# Household-Specific Regressions Using Clickstream Data

**Avi Goldfarb and Qiang Lu**

*Abstract.* This paper makes three contributions: (1) the paper provides a better understanding of online behavior by showing the main drivers of Internet portal choice, (2) the rich data allow for a deeper understanding of brand substitution patterns than previously possible and (3) the paper introduces a wider statistics community to a new data opportunity and a recently developed method.

*Key words and phrases:* Internet portals, clickstream data, online behavior.

## 1. INTRODUCTION

The rapid growth of the Internet since the late 1990s has been an amazing opportunity for marketing and statistics research. We can now observe consumers in a new kind of environment where information search is relatively inexpensive. There are both direct and indirect reasons to be interested in online behavior. The direct reason is that the Internet has become an important element of daily life for millions of people. A growing body of evidence shows that online behavior differs from offline behavior (e.g., Park and Fader, 2004; Chen and Hitt, 2002; Danaher, Wilson and Davis, 2003). Therefore, it is essential for companies with an online presence to understand what drives online choices. The indirect reason is that Internet data and detail on online choices can inform researchers about more general economic questions such as the cost of search. [Ellison and Ellison (2005) provide a review of how Internet research has informed a number of general economics and marketing questions.]

Studies of online behavior have one important advantage over previous studies of behavior: rich, detailed data. Clickstream data tracks the website visit history for each online user, thereby enabling a detailed analysis of user choice. One advantage of clickstream data sets is that they are much larger than the data sets typically used to examine consumer behavior. A typical household may buy consumer goods like ketchup and detergent every month or two; however, this same household may visit a given website ten times each week or more. Clickstream data can provide us much more information about Internet users, and can therefore allow us to do more complicated analysis than with traditional data. At the same time, this rich data means it is difficult for managers to use standard methods to quickly arrive at useful insights using their desktop computers.

This paper uses detailed clickstream data on Internet portals to better understand online choices. [In this study, we define an Internet portal as any site that classifies content and primarily presents itself as a one-stop point-of-entry to content on the web (Hargittai, 2000).] The paper makes three contributions. First, it provides a better understanding of online behavior by showing the main drivers of Internet portal choice. This relates to the direct reason to study the Internet. In this paper, we show that the success of previous searches is a particularly important driver of website choice. The ability to provide deep searches with many results is less important. This may partially explain how Google, with its page-rank technology and "I feel lucky" button, grew to be dominant in search despite being a late entrant to the portal/search market.

*Avi Goldfarb is Assistant Professor, Rotman School of Management, University of Toronto, Toronto, Ontario, Canada MSS 3E6 e-mail: agoldfarb@rotman.utoronto.ca. Qiang (Steven) Lu is Assistant Professor, School of Business, University of Sydney, New South Wales 2006, Australia e-mail: s.lu@econ.usyd.edu.au.*

Second, the rich data allow for a deeper understanding of brand substitution patterns than previously possible. This relates to the indirect reason to study the Internet. In our analysis, we show that consumers who rarely change portals often prefer either Yahoo!, the most popular portal, or Go2net and Ask Jeeves, both among the least popular portals. In contrast, those that prefer other popular portals such as Excite, MSN, and especially AOL are more likely to change the portal they visit. Many brands seem to have less loyal consumers. Both product features and habit formation are known to be important parts of brand building (e.g., Aaker, 1991). Clickstream data provides a method for isolating the relationship between habit formation and brand strength that was not previously possible. In this context, habit formation is shown to be relatively uncorrelated with brand preferences.

Third, the paper introduces a wider statistics community to a new data opportunity and a recently developed method. We examine the choices of Internet portals by 2517 online households over three months. The average household makes 312 portal visits in this time. (The data are available to university-affiliated researchers. If interested in using this data, email Avi Goldfarb at agoldfarb@rotman.utoronto.ca.) This rich data set means that we can conduct household-specific regressions. By "household-specific regressions" we mean that we run a separate conditional logit regression on the portal choice of each household in the data. Following Goldfarb (2006), this method allows for more flexible (semiparametric) substitution patterns than the panel methods typically used to study offline consumer choice because all coefficients are allowed to vary at the household level without any assumed distribution. The method has the added advantage of being computationally inexpensive and relatively easy to understand. This means that managers can use it to quickly understand the behavior of their customers.

The next section details the data set used in the analysis. This is followed by a description of the methodology, the results and a brief conclusion.

## 2. DATA

The data set, provided by Plurimus Corporation, consists of 2,645,778 website visits by 2651 households from December 27, 1999 to March 31, 2000. For our analysis, we have chosen the top 15 Internet portals visited by 2517 different users with 784,882 observations. This data set is larger than those typically used in marketing and economics to analyze consumer choice. The raw data show the household identification, the website visited, the time of day the visit started and ended (in seconds) and the number of pages viewed at the website.

Despite its richness, the data set has four main limitations. First, it is not geographically representative. Most respondents are drawn from Pennsylvania, Texas, Florida and North Carolina. Second, it contains few users at-work. Third, the data are collected at the household level. And fourth, the data do not contain information on households the first time they go online. These limitations mean that results should be extended to different geographic distributions and at-work users with caution. Furthermore, the third limitation means that consumer loyalty may be underestimated if different household members have different preferences. The third and fourth limitations together mean that it is difficult to study consumer learning.

From this data, measures of loyalty, past number of pages viewed at the website and past search failure are constructed. Loyalty is measured as a dummy variable for whether that website was visited on the previous choice occasion: it is equal to 1 if the website was visited on the previous occasion and zero otherwise. [In the marketing literature, the variable measuring the previous choice is typically called "loyalty" (e.g., Guadagni and Little, 1983). More formally, this variable measures "true state dependence."] Past number of pages viewed is measured by the number of pages viewed at the website on the previous visit. Past search failure is measured by whether the previous search is repeated. In particular, if a household conducts two searches within five minutes for the same goal (over 100 possible goals were identified in the raw data) then it is considered a repeated search. Furthermore, if a household visited two portal sites in a row, and there were less than five minutes between visits, then the first search is considered to have been repeated. It is pages viewed and repeated search during the previous visit to a portal that are relevant for the current choice of which website to visit. Therefore we include data on the experience of the household during its previous visit to the website in the analysis. The use of data on the previous experience at the website means that a correction is required for households that have not yet visited a given website. A *Missing Data* variable is used for this correction. It is equal to 1 if there is no information about a household's previous visit to a website (because the household has not yet visited it) and zero otherwise. This variable serves as a control and has no clear economic interpretation. Table 1 provides univariate descriptive statistics. Table 2 provides pairwise correlation coefficients based on household-level aggregates. In particular, the average values for the variables were calculated at the household level. These are the correlations of the values across the 2517 households in the sample.

Table 1
*Basic descriptive statistics*

| Variable | Mean | Std. dev. | Min. | Max. |
|---|---|---|---|---|
| Loyalty | 0.10 | 0.31 | 0 | 1 |
| Last Search Repeated | 0.25 | 0.43 | 0 | 1 |
| Ln(Last # Pages) | 0.59 | 0.86 | 0 | 7.46 |
| Missing Data | 0.23 | 0.42 | 0 | 1 |

## 3. METHOD

The decision to visit an Internet portal is modeled as a standard economic discrete-choice problem. Users visit the portal that they expect to give them the highest utility. The utility to household $i$ of visiting portal $j$ at time $t$ is

$$u_{ijt} = X_{ijt}\beta_i + \varepsilon_{ijt}, \tag{1}$$

where $X_{ijt}$ includes loyalty, whether the last search was repeated, the number of pages viewed during the last visit to the website (naturally logged), a missing data measure and the brand dummy variables; $\beta_i$ is a vector of household-specific coefficients; and $\varepsilon_{ijt}$ is an idiosyncratic error term. Assuming that $\varepsilon_{ijt}$ follows a type-II extreme value distribution over $j$ and $t$, we can estimate equation (1) using a separate conditional logit regression (McFadden, 1974) for each household in the data. Therefore the probability that household $i$ visits portal $j$ at time $t$ is modeled as

$$P_{ijt} = \frac{\exp(X_{ijt}\beta_i)}{\sum_{k=1}^{J_i}\exp(X_{ikt}\beta_i)}. \tag{2}$$

Typically, the parameters $\beta_i$ are fixed across households or assumed to follow a known (typically normal) distribution. Household-specific regressions allow the coefficients to vary nonparametrically across households and are computationally much less intensive. [A number of previous studies have recommended running regressions on the time-varying dimension of a panel data set when there is sufficient data. Fischer and Nagin (1981) conducted experiments that showed coefficients vary across individuals. Pesaran and Smith (1995) examined employment in 38 industries and concluded (page 102) that the "lesson for applied work is that when large T panels are available, the individual micro-relations should be estimated separately." Elrod and Häubl (1998), however, highlight the shortcomings of household-specific regressions: inefficiency, overestimation of population variance and an inability to predict out-of-sample. Therefore, household-specific regressions are best used to understand the underlying drivers of behavior in very large data sets.] In particular, a separate vector of coefficients is estimated for each household. [Numerous previous studies have shown the logit model to fit consumer choice behavior well (e.g., Guadagni and Little, 1983; Jain, Vilcassim and Chintagunta, 1994). Still, these studies do not examine the fit at that household level. Goldfarb (2006) shows that taken together, household-specific regressions provide a better fit than aggregated methods. It is, however, likely that the fit is poor for some specific households in the data set.] Furthermore, the discrete-choice random coefficients models used in these settings may take weeks to converge with millions of observations (and only if the researchers have sufficient RAM in their computers). This matters if managers plan to use clickstream data analysis to make business decisions. Goldfarb (2006) shows that the coefficient estimates are consistent and directly comparable under a mild set of assumptions. Household-specific regressions require one key additional restriction because most households do not visit all 15 portals: only those websites that a particular household actually visits over the course of the sample are included in the analysis. For example, if household $i$ only visits Yahoo!, MSN, Excite and Go.com, then only these four portals are included in the regressions for that household. Implicitly, the coefficients on the dummies for the other portals approach negative infinity.

The likelihood function for household $i$ is

$$\prod_{t=1}^{T_i}\prod_{j=1}^{J_i}(P_{ijt})^{d_{ijt}}, \tag{3}$$

Table 2: Correlation coefficients of household-level data

| | **Total number pages visited** | **Average number pages visited** | **% searches repeated** | **% Yahoo** | **% MSN** | **% Netscape** | **% Excite** | **% AOL** | **% Altavista** | **% Iwon** | **% Lycos** | **% MyWay** | **% Go** | **% Hotbot** | **% Snap** | **% Go2Net** | **% Goto** | **% Ask Jeeves** |
|---|---|---|---|---|---|---|---|---|---|---|---|---|---|---|---|---|---|---|
| Total # pages | 1 | | | | | | | | | | | | | | | | | |
| Average # pages | 0.01 | 1 | | | | | | | | | | | | | | | | |
| Search repeated | 0.25[a] | −0.04[a] | 1 | | | | | | | | | | | | | | | |
| %Yahoo | 0.06[a] | −0.23[a] | −0.09[a] | 1 | | | | | | | | | | | | | | |
| % MSN | −0.06[a] | −0.05[a] | −0.18[a] | −0.37[a] | 1 | | | | | | | | | | | | | |
| % Netscape | −0.01 | 0.12[a] | 0.01 | −0.27[a] | −0.28 | 1 | | | | | | | | | | | | |
| % Excite | 0.07[a] | 0.07[a] | 0.13[a] | −0.16[a] | −0.14[a] | −0.06[a] | 1 | | | | | | | | | | | |
| % AOL | −0.02 | 0.12[a] | −0.03 | −0.15[a] | −0.10[a] | −0.11[a] | −0.04 | 1 | | | | | | | | | | |
| % Altavista | −0.02 | 0.13[a] | 0.01 | −0.14[a] | −0.09[a] | −0.12[a] | −0.06[a] | −0.05[a] | 1 | | | | | | | | | |
| % Iwon | 0.11[a] | 0.16[a] | 0.01 | −0.10[a] | −0.08[a] | −0.06[a] | −0.03 | −0.05[a] | −0.03 | 1 | | | | | | | | |
| % Lycos | −0.04[a] | −0.07[a] | 0.21[a] | −0.14[a] | −0.10[a] | −0.01 | −0.03 | −0.05[a] | −0.04[a] | −0.01 | 1 | | | | | | | |
| % MyWay | −0.01 | 0.03 | 0.00 | −0.13[a] | −0.08[a] | −0.08[a] | −0.05[a] | −0.05[a] | −0.05[a] | −0.03 | −0.03 | 1 | | | | | | |
| % Go | −0.07[a] | 0.09[a] | −0.00 | −0.11[a] | −0.05[a] | −0.10[a] | −0.05[a] | −0.03 | −0.01 | −0.04 | −0.02 | −0.04 | 1 | | | | | |
| % Hotbot | −0.03 | −0.04[a] | 0.18[a] | −0.13[a] | −0.12[a] | 0.01 | −0.02 | −0.04[a] | −0.04[a] | −0.03 | 0.35[a] | −0.03 | −0.04[a] | 1 | | | | |
| % Snap | 0.06[a] | 0.01 | 0.09[a] | −0.09[a] | −0.07[a] | −0.02 | −0.09 | −0.03 | −0.01 | −0.01 | −0.02 | −0.02 | −0.03 | −0.02 | 1 | | | |
| % Go2Net | −0.02 | −0.08[a] | 0.12[a] | −0.10[a] | −0.07[a] | −0.04[a] | −0.13 | −0.02 | −0.04 | −0.03 | −0.03 | −0.04 | −0.04 | −0.03 | −0.02 | 1 | | |
| % Goto | −0.06[a] | 0.00 | 0.21[a] | −0.14[a] | −0.01 | −0.04[a] | −0.03 | −0.03 | −0.03 | −0.03 | 0.02 | −0.01 | 0.03 | −0.01 | −0.03 | 0.22[a] | 1 | |
| % Ask Jeeves | −0.03 | 0.02 | 0.02 | −0.06[a] | −0.08[a] | −0.04[a] | −0.03 | −0.01 | −0.00 | −0.02 | −0.02 | −0.02 | 0.02 | −0.02 | −0.01 | −0.02 | −0.02 | 1 |

[a]Significant at 95% confidence level.

where $d_{ijt}$ is equal to 1 if household $i$ purchased brand $j$ at time $t$, 0 otherwise. Maximum likelihood estimation is applied for each household. Thus, for 2517 households, 2517 sets of parameters are estimated.

## 4. RESULTS

The results are presented in Tables 3 and 4 and in Figures 1 to 4. The regressions led to 2517 different vectors of parameter estimates. Table 3 presents a univariate description of these parameter estimates. It shows the mean of each coefficient across households, the standard error of the mean, the standard deviation of the coefficients and the percentage of coefficients that are significantly different from zero with 95% confidence (positive and negative). Table 3 shows that loyalty and repeated search are important factors in portal choice. The result on repeated search, relative to pages viewed, is particularly interesting. It suggests that it is much more important for portals to direct people to the right website than for portals to provide many results. Information quality is much more important than information quantity. This may partially explain the rise of Google in the years since the data was collected. (Google is the 17th most popular portal in the data and is therefore not included in the study.) Even before Google had a large database, its technology was particularly good at ordering results. This likely reduced the frequency of repeated searches. The coefficients on the portal dummies are all presented relative to Yahoo!. [The portal dummy coefficient distributions are based only on those 2206 households that visited Yahoo! at least once. This ensured that household-level coefficients would be comparable since they would all have the same base (Yahoo!).] Since the mean coefficients are all negative, this suggests that Yahoo! is, on average, the preferred portal by a substantial margin (controlling for loyalty, search repeats and pages viewed). MSN and Netscape also have substantially higher coefficients. The remaining twelve portals have relatively similar coefficient values. This characterization is subject to the important caveat that there is considerable heterogeneity across responses. The standard deviations of the coefficients are high relative to the means.

Table 4 presents the pairwise correlation coefficients of the parameter estimates. It shows a number of interesting results. First, there is a strong negative correlation between Last Search Repeated and Ln(Last # Pages), which is quite intuitive: the more a user cares about accuracy of the search results, the less he cares about the number of search results. Figure 1 presents a scatterplot of these variables. Second, Loyalty does not seem to be correlated with Last Search Repeated or with Ln(Last # Pages). Consumers who tend to return to the same portals are no more or less likely to care about search accuracy or search depth.

Third, Loyalty has an interesting relationship with preferences for the various brands. Consumers who rarely change portals have lower opinions of MSN, Excite, Altavista, and especially AOL relative to their opinions of Yahoo!. Loyal households prefer Yahoo! and Netscape relative to the other major brands. Figure 2 displays a scatterplot of the Loyalty coefficients relative to the coefficients of four popular portals. Interestingly, loyal households are also particularly likely to visit unpopular portals such as Go2Net, Hotbot and Ask Jeeves. Figure 3 displays the scatterplot. Many strong brands seem to have less loyal consumers, suggesting that in this context brand strength relates to product features more than simply habit formation. While the extant literature acknowledges the roles of both habit formation and product features in building brands (e.g., Aaker, 1991), the rich clickstream data set used here allows for an assessment of the correlation between brand strength and loyalty that was not previously possible.

Fourth, brand preferences relative to Yahoo! are positively correlated. This result emphasizes Yahoo!'s dominance and is shown in a scatterplot in Figure 4. Unsurprisingly, many brands that are owned by the same company or are linked to each other have highly correlated preference coefficients. For example, preferences for Lycos are highly correlated with Hotbot. Both had the same owner at the time. More interesting is a comparison of the MSN and Netscape correlation coefficients. Those brand coefficients that are highly correlated with Netscape are less correlated with MSN and vice versa. Overall, these results suggest that some portals are substitutes for each other in that they may be used for the same purposes. On the other hand, some portals are complements for each other. People may visit Lycos for searches about celebrities and then Ask Jeeves for science-related questions.

Last, there is a strong negative correlation between Last Search Repeated and Netscape, which

TABLE 3

*Descriptive statistics of the 2517 household-specific coefficient vectors*

| **Variable** | **Mean coefficient** | **Standard error of the mean** | **Standard deviation of the coefficients** | **% Significantly positive (95%)** | **% Significantly negative (95%)** |
|---|---|---|---|---|---|
| Variables | | | | | |
| Loyalty | 0.98 | 0.01 | 2.78 | 69.62 | 4.41 |
| Last Search Repeated | −0.39 | 0.01 | 1.48 | 3.19 | 37.14 |
| Ln(Last # Pages) | 0.03 | 0.01 | 2.60 | 12.58 | 6.73 |
| Missing Data | −0.63 | 0.02 | 5.07 | 7.35 | 21.96 |
| Brand dummies | | | | | |
| MSN | −0.81 | 0.01 | 4.89 | 28.60 | 41.77 |
| Netscape | −1.57 | 0.01 | 5.25 | 23.53 | 48.08 |
| Excite | −2.12 | 0.01 | 3.08 | 8.30 | 59.43 |
| AOL | −2.07 | 0.01 | 2.47 | 14.40 | 61.80 |
| Altavista | −2.60 | 0.01 | 2.34 | 4.95 | 64.77 |
| Iwon | −3.00 | 0.01 | 2.58 | 5.13 | 69.25 |
| Lycos | −1.96 | 0.01 | 4.23 | 10.66 | 57.78 |
| MyWay | −2.38 | 0.01 | 3.63 | 6.90 | 62.99 |
| Go | −2.80 | 0.01 | 3.49 | 7.26 | 62.38 |
| Hotbot | −2.28 | 0.01 | 2.18 | 7.07 | 61.41 |
| Snap | −2.89 | 0.01 | 3.23 | 6.00 | 66.17 |
| Go2Net | −2.14 | 0.01 | 3.10 | 10.52 | 58.08 |
| Goto | −2.69 | 0.01 | 2.16 | 11.96 | 64.74 |
| Ask Jeeves | −2.71 | 0.01 | 3.32 | 17.63 | 57.93 |

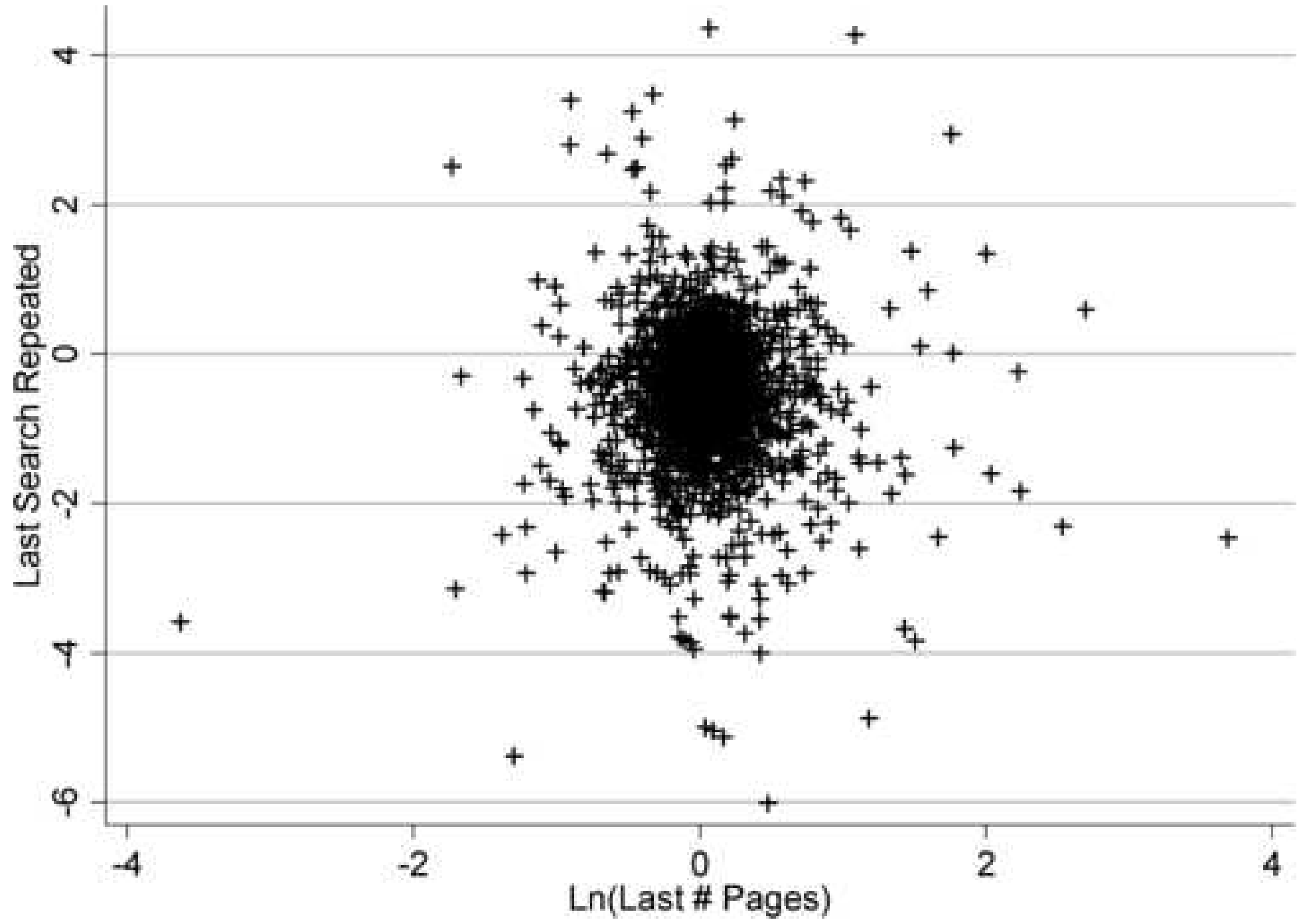

FIG. 1. *Scatterplot of the coefficients Ln(Last # Pages) and Last Search Repeated.*

Table 4: Correlation coefficients of the household-specific coefficients

| | Loyalty | Last Search Repeated | Last # Pages visited | Missing | MSN | Netscape | Excite | AOL | Altavista | Iwon | Lycos | MyWay | Go | Hotbot | Snap | Go2Net | Goto | Ask Jeeves |
|---|---|---|---|---|---|---|---|---|---|---|---|---|---|---|---|---|---|---|
| Total # pages | 1 | | | | | | | | | | | | | | | | | |
| Last Search Repeated | −0.01 | 1 | | | | | | | | | | | | | | | | |
| Last # pages visited | −0.02 | **−0.90**[a] | 1 | | | | | | | | | | | | | | | |
| Missing | 0.33[a] | −0.73[a] | 0.80[a] | 1 | | | | | | | | | | | | | | |
| MSN | **−0.08**[a] | −0.15[a] | 0.21[a] | 0.02 | 1 | | | | | | | | | | | | | |
| Netscape | **0.02** | **−0.99**[a] | **0.97**[a] | 0.96[a] | 0.61[a] | 1 | | | | | | | | | | | | |
| Excite | **−0.16**[a] | −0.53[a] | 0.35[a] | 0.28[a] | 0.73[a] | 0.88[a] | 1 | | | | | | | | | | | |
| AOL | **−0.77**[a] | −0.39[a] | −0.79[a] | −0.81[a] | 0.67[a] | 0.36[a] | 0.94[a] | 1 | | | | | | | | | | |
| Altavista | **−0.16**[a] | −0.43[a] | −0.03 | 0.04 | 0.70[a] | 0.57[a] | −0.20[a] | 0.56[a] | 1 | | | | | | | | | |
| Iwon | **−0.17**[a] | −0.81[a] | 0.34[a] | 0.80[a] | 0.91[a] | 0.67[a] | 0.60[a] | 0.49[a] | 0.89[a] | 1 | | | | | | | | |
| Lycos | −0.46[a] | −0.30[a] | 0.30[a] | 0.22[a] | 0.86[a] | 0.78[a] | 0.83[a] | 0.47[a] | 0.58[a] | 0.70[a] | 1 | | | | | | | |
| MyWay | 0.28[a] | −0.47[a] | −0.29[a] | 0.23[a] | 0.45[a] | 0.85[a] | 0.45[a] | 0.45[a] | 0.73[a] | 0.72[a] | 0.59[a] | 1 | | | | | | |
| Go | −0.18[a] | 0.17[a] | −0.26[a] | −0.51[a] | 0.74[a] | 0.87[a] | 0.72[a] | 0.29[a] | 0.79[a] | 0.58[a] | 0.77[a] | 0.52[a] | 1 | | | | | |
| Hotbot | −0.82[a] | −0.40[a] | −0.55[a] | −0.82[a] | 0.82[a] | 0.51[a] | 0.86[a] | 0.49[a] | 0.62[a] | 0.65[a] | **0.88**[a] | 0.78[a] | 0.85[a] | 1 | | | | |
| Snap | −0.12[a] | −0.52[a] | −0.45[a] | −0.19[a] | 0.26[a] | 0.57[a] | 0.57[a] | 0.51[a] | 0.54[a] | 0.60[a] | 0.83[a] | 0.56[a] | 0.46[a] | 0.83[a] | 1 | | | |
| Go2Net | **0.40**[a] | −0.05 | −0.82[a] | 0.17[a] | 0.55[a] | 0.81[a] | 0.61[a] | 0.93[a] | 0.67[a] | 0.52[a] | 0.74[a] | 0.48[a] | 0.44[a] | 0.57[a] | 0.65[a] | 1 | | |
| Goto | 0.36[a] | −0.41[a] | 0.56[a] | 0.81[a] | 0.75[a] | 0.85[a] | 0.40[a] | 0.49[a] | 0.81[a] | 0.52[a] | 0.37[a] | 0.61[a] | 0.53[a] | 0.51[a] | 0.35[a] | 0.42[a] | 1 | |
| Ask Jeeves | **0.43**[a] | 0.36[a] | 0.06 | 0.50[a] | 0.58[a] | 0.49[a] | 0.42[a] | 0.48[a] | 0.33[a] | 0.51[a] | 0.43[a] | 0.49[a] | 0.67[a] | 0.59[a] | 0.38[a] | 0.48[a] | 0.47[a] | 1 |

[a]Correlations mentioned in the text are written in **bold**. Significant at 95% confidence level.

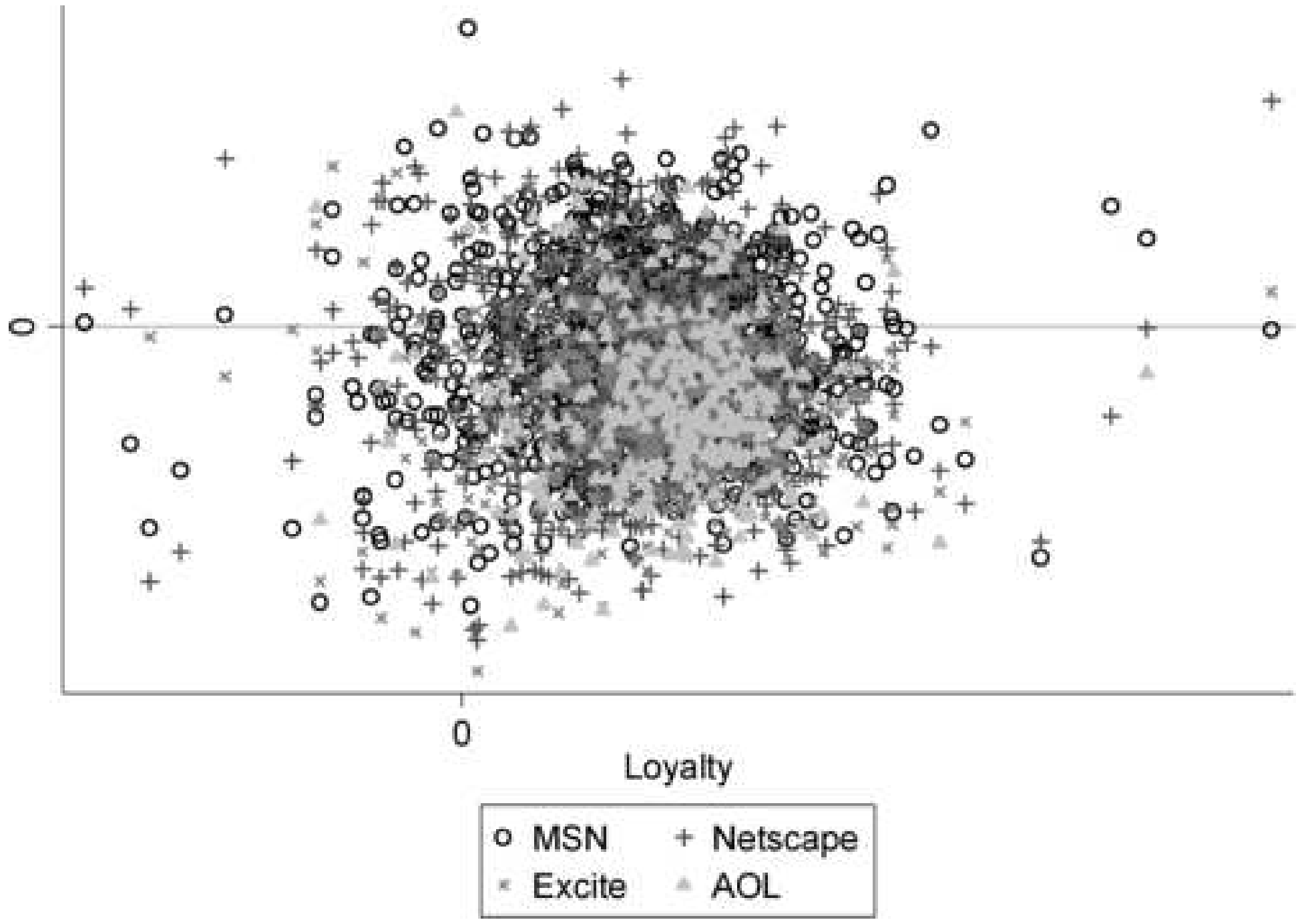

Fig. 2. *Scatterplot of the Loyalty coefficients relative to popular portal dummies.*

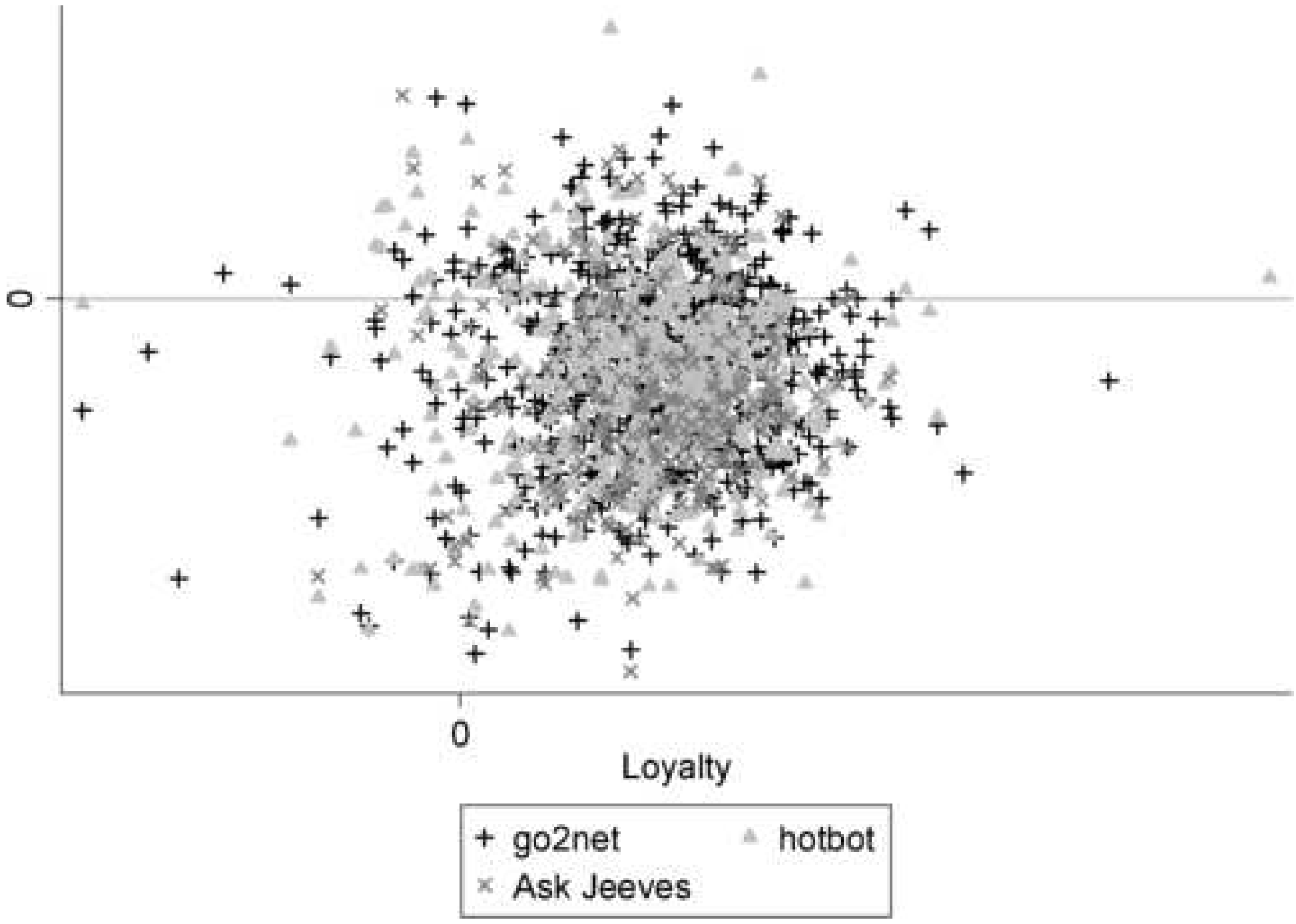

Fig. 3. *Scatterplot of the Loyalty coefficients relative to the less popular portal dummies.*

means the more important the information accuracy to a household, the less the household visits Netscape; and there is a strong positive correlation between Ln(Last # Pages) and Netscape, which means the more a household cares about the number of search results, the better image it has of Netscape. At the time, Netscape ordered its search results by the search technology used rather than by relevance.

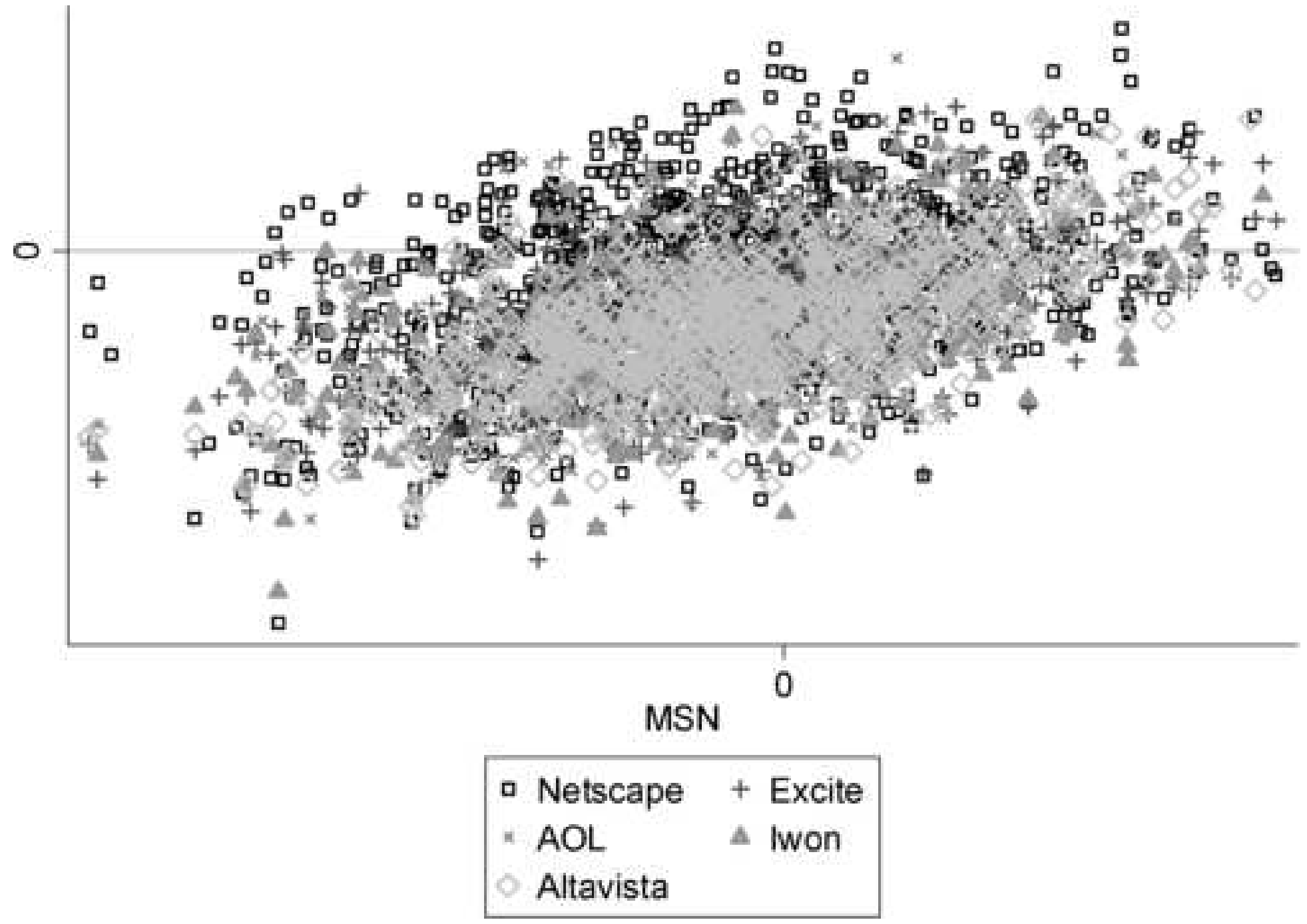

Fig. 4. *Scatterplot of the coefficients of the MSN dummy relative to five other portal dummies.*

This provides an important reality check on the results suggesting people who preferred many results to targeted results preferred a portal that provided many untargeted results.

In summary, examining the relationship between the parameter estimates across households provides a rich insight into consumer behavior that was not previously possible.

## 5. CONCLUSIONS

In summary, this paper provides a better understanding of Internet portal choice. It also shows that the household-level correlation of brand preferences and loyalty provides insights into brand building more generally. Finally, it introduces a wider audience to a rich clickstream data set and to household-specific regression. Clickstream data is an exciting opportunity for statisticians. Future research can proceed in a number of different directions. It can explore how behavior has changed as people become more comfortable online. It can explore how the online environment affects purchase behavior. It can use the data to develop new statistical tools for large data sets. And it can better inform our understanding of general economic and psychological issues. Going forward, there are a number of exciting research opportunities arising from the availability of clickstream data.

## ACKNOWLEDGMENTS

We thank the Guest Editors Wolfgang Jank and Galit Shmueli, the referees and participants at the First StatsChallenges in Ecommerce conference for helpful comments. Jin Wang provided excellent research assistance. Funding was provided by SSHRC Grant 538-02-1013. Please address correspondence to agoldfarb@rotman.utoronto.ca. All opinions and errors are our own.